%% Plain TeX
\magnification=1200
\font\titlefont=cmcsc10 at 12pt
\hyphenation{moduli}
 
%
%     *****Part of MSSYMB.Tex*****
%
\catcode`\@=11
\font\tenmsa=msam10
\font\sevenmsa=msam7
\font\fivemsa=msam5
\font\tenmsb=msbm10
\font\sevenmsb=msbm7
\font\fivemsb=msbm5
\newfam\msafam
\newfam\msbfam
\textfont\msafam=\tenmsa  \scriptfont\msafam=\sevenmsa
  \scriptscriptfont\msafam=\fivemsa
\textfont\msbfam=\tenmsb  \scriptfont\msbfam=\sevenmsb
  \scriptscriptfont\msbfam=\fivemsb
\def\hexnumber@#1{\ifcase#1 0\or1\or2\or3\or4\or5\or6\or7\or8\or9\or
      A\or B\or C\or D\or E\or F\fi }
 
%  The following 13 lines establish the use of the Euler Fraktur font.
\font\teneuf=eufm10
\font\seveneuf=eufm7
\font\fiveeuf=eufm5
\newfam\euffam
\textfont\euffam=\teneuf
\scriptfont\euffam=\seveneuf
\scriptscriptfont\euffam=\fiveeuf
% Use the next 4 lines with AMS-TeX:
%\def\frak{\relaxnext@\ifmmode\let\next\frak@\else
% \def\next{\Err@{Use \string\frak\space only in math mode}}\fi\next}
%\def\goth{\relaxnext@\ifmmode\let\next\frak@\else
% \def\next{\Err@{Use \string\goth\space only in math mode}}\fi\next}
% Use the next 4 lines if NOT using AMS-TeX:
\def\frak{\ifmmode\let\next\frak@\else
 \def\next{\errmessage{Use \string\frak\space only in math mode}}\fi\next}
\def\goth{\ifmmode\let\next\frak@\else
 \def\next{\errmessage{Use \string\goth\space only in math mode}}\fi\next}
\def\frak@#1{{\frak@@{#1}}}
\def\frak@@#1{\fam\euffam#1}
%  End definition of Euler Fraktur font.
 
\edef\msa@{\hexnumber@\msafam}
\edef\msb@{\hexnumber@\msbfam}
\mathchardef\square="0\msa@03
\mathchardef\subsetneq="3\msb@28
\mathchardef\ltimes="2\msb@6E
\mathchardef\rtimes="2\msb@6F
\def\Bbb{\ifmmode\let\next\Bbb@\else
\def\next{\errmessage{Use \string\Bbb\space only in math mode}}\fi\next}
\def\Bbb@#1{{\Bbb@@{#1}}}
\def\Bbb@@#1{\fam\msbfam#1}
\catcode`\@=12
%
%     *****End of MSSYMB.Tex*****
%

\def\E{{\Bbb E}}

\def\Q{{\Bbb Q}}

\def\pn{\par\noindent}
\def\no{\noindent}
\def\mbar#1{{\overline{\cal M}}_{#1}}
\def\Mg{{\cal M}_g}
\def\la#1{\lambda_{#1}}
\def\ka#1{\kappa_{#1}}

\def\m#1{{\cal M}_{#1}}

\def\a#1{{\cal A}_{#1}}
\def\j#1{{\cal J}_{#1}}
\def\ch#1{ch_{#1}}
\def\ll{\langle}
\def\rr{\rangle}

\def\Dirr{\Delta_{\rm irr}}
\def\dirr{\delta_{\rm irr}}
\def\nbar{\underline{n}}
 
\vskip 6.5pc
\noindent
\font\eighteenbf=cmbx10 scaled\magstep2
\vskip 2.0pc
\centerline{\eighteenbf Algorithms for computing intersection numbers }
\medskip
\centerline{\eighteenbf on moduli spaces
of curves, with an application	}
\medskip
\centerline{\eighteenbf to the class of the locus of Jacobians	}
\noindent
\vskip 2pc
\font\titlefont=cmcsc10 at 11pt
\centerline{\titlefont Carel Faber}
\vskip 2.0pc

\bigskip
The purpose of this note is first of all
to explain how intersection numbers of divisors
on the moduli spaces $\mbar{g,n}$ of stable pointed curves can be computed.
The Witten conjecture, proven by Kontsevich, gives a recipe to compute the
intersection numbers of the $n$ basic line bundles on $\mbar{g,n}\,$. As we
will see, the knowledge of these numbers allows one to compute all other
intersection numbers of divisors as well. That this is possible was pointed
out to me by Rahul Pandharipande. Earlier, Eduard Looijenga had made a remark
that went a long way in the same direction.

After describing the various divisors on $\mbar{g,n}\,$, we proceed to
discuss the algorithm for the computation of their intersection numbers.
We have implemented the algorithm and we will discuss this implementation
as well as the results obtained this way. E.g., we computed all intersection
numbers on $\mbar{g}$ for $g\le6$. (Copies of the program and some data
computed with it are available from the author.)

A refined version of the algorithm requires to take certain 
higher codimensional classes, introduced by Mumford and Arbarello-Cornalba,
into account; it computes all intersection numbers of these classes and
divisors.
Recently, we realized that the Chern classes of the Hodge bundle can be
taken along as well; hence all intersection numbers of Mumford's
tautological classes and divisors can be computed. This has several
applications. In \S4 we discuss one application in detail: 
the calculation of the class of the locus of Jacobians in
the moduli space of principally polarized abelian varieties of dimension $g$
(projected in the tautological ring).
This class was classically known for $g=4$ and was computed by us with
ad hoc methods for $g=5\,$; the new method allows in principle to compute
the class for all $g$, in practice currently for $g\le7$.
Other applications may be found in the recent papers of Graber and
Pandharipande [GP] and Kontsevich and Manin [KM].

\bigskip
\no {\bf \S1. Line bundles and divisors on $\mbar{g,n}\,$.}

\medskip
For non-negative integers $g$ and $n$ with $2g-2+n>0$, denote by
$\mbar{g,n}$ the moduli space of stable $n$-pointed curves of genus $g$,
over an algebraically closed field $k$. This is the Deligne-Mumford
compactification of the moduli space $\m{g,n}$ of smooth $n$-pointed curves
$(C;x_1,\dots,x_n)$ of genus $g$ (with $x_i\neq x_j$ if $i\neq j$).
We will consider certain classes in the rational Picard group of
$\mbar{g,n}\,$. First, for $1\le i\le n$, let $\psi_i$ denote the first
Chern class of the line bundle whose fiber at a stable $n$-pointed curve
$(C;x_1,\dots,x_n)$ is the cotangent space to $C$ at $x_i\,$; i.e.,
$\psi_i=c_1(\sigma_i^*(\omega_{\pi_{n+1}}))$, where $\pi_{n+1}:
\mbar{g,n+1}\to\mbar{g,n}$ is the morphism obtained by forgetting the $(n+1)$-st
marked point (the universal curve, cf.~[Kn~1]), $\omega_{\pi_{n+1}}$
is the relative dualizing sheaf, and $\sigma_1,\dots,\sigma_n$ are the
natural sections of $\pi_{n+1}$ (the image of a stable $n$-pointed curve
under $\sigma_i$ is the stable $(n+1)$-pointed curve obtained by attaching
a 3-pointed rational curve at the $i$-th point and considering the remaining
2 points on that curve as the $i$-th and $(n+1)$-st point).
Next, following [AC~2], \S1, we define $\ka1=\pi_{n+1,*}(K^2)$, with
$K=c_1(\omega_{\pi_{n+1}}(\sum_{i=1}^n D_i))$, where $D_i$ is the divisor
that is the image of the section $\sigma_i\,$, for $1\le i\le n$.
Note that it is a consequence of results of Harer (cf.~[Ha~1], [Ha~2],
[AC~1]) that over ${\Bbb C}$ the restrictions to $\m{g,n}$
of the classes $\ka1$ and $\psi_1,\dots,\psi_n$ generate the rational Picard
group of $\m{g,n}\,$. 

To get generators for the rational Picard group of $\mbar{g,n}\,$, we have
to add the fundamental classes of the boundary divisors.
Exactly when $g>0$, there is a boundary component whose generic point
corresponds to an irreducible singular curve. It is the image of
$\mbar{g-1,n+2}$ under the degree $2$ map that identifies the $(n+1)$-st
and $(n+2)$-nd point on each curve. Following [AC~2], we denote
this locus by $\Dirr$ and its class in the Picard group by
$\dirr\,$. (For $g=0$ this class is $0$ by definition.)

The other boundary components parametrize reducible singular
curves. The generic point of such a component corresponds to a curve
with two irreducible components $C_1$ and $C_2\,$, 
with genera $g_1$ and $g_2$ satisfying
$g_1+g_2=g$, and labelled by subsets $N_1$ and $N_2$ of
$\nbar=\{1,2,\dots,n\}$ satisfying $N_1\coprod N_2=\nbar$
that correspond to the marked points on the two components. 
All partitions of $g$ and $\nbar$ that lead to a stable curve
occur; this just translates as the condition $|N_i|\ge2$ when
$g_i=0$. Such a boundary component is the image of
$\mbar{g_1,|N_1|+1}\times\mbar{g_2,|N_2|+1}$ under the
natural map that identifies the two `extra' points and
labels the $|N_i|$ remaining points on $C_i$ with the labels from $N_i\,$.
We have chosen to denote this boundary component in the case $n>0$
by $\Delta_{g_i,N_i}$ where $N_i$ is the subset of $\nbar$ containing
$1$. 
In the case $n=0$ the $N_i$ are empty and we may drop them in the notation;
note that $\Delta_{g_1}=\Delta_{g_2}$ and that this component is usually
denoted as $\Delta_{\min(g_1,g_2)}\,$.

Although this will play no role in the sequel, we point out that the
classes in the rational Picard group that we have introduced so far, are
independent whenever $g\ge3$. For $g=2$, there is one relation, originating
from the fact that $\kappa_1$ on $\mbar{2}$ comes from the boundary.
For $g=1$, both $\kappa_1$ and the $\psi_i$'s come from the boundary;
the boundary components are independent. For $g=0$, the
boundary components generate, but are not independent; the relations
arise from the various projections to $\mbar{0,4}$ and the equivalence
of its 3 boundary cycles (cf.~[Ke]). Let me add that although the
results mentioned in this paragraph are no doubt correct, it appears that
for some of the statements, references containing proofs are not (yet)
available.

There is one other divisor class which will be most useful: $\la1\,$, the
first Chern class of the Hodge bundle. The Hodge bundle is the locally
free rank $g$ sheaf (on the moduli functor) whose fiber at a curve
$C$ is $H^0(C,\omega_C)$. So it is $0$ in genus 0, while it is a pull-back
from $\mbar{1,1}$ resp.~$\mbar{g}$ in case $g=1$ resp.~$g\ge2$.

\bigskip
\no {\bf \S2. The idea of the algorithm.}

\medskip
Suppose given a monomial of degree $3g-3+n$ in the divisor classes
$\ka1$, $\psi_1$, \dots, $\psi_n$, $\dirr$ and the $\delta_{g_i,N_i}$ on
$\mbar{g,n}\,$; we want to compute the corresponding intersection number.
We will interpret the divisor classes as classes on the moduli functor. 
In the case of a boundary divisor, this means that we divide the
usual fundamental class by the order of the automorphism group of the
generic curve parametrized by the divisor. These divisor classes will be
denoted by $\delta_{\dots}\,$, to distinguish them from the actual
boundary divisors $\Delta_{\dots}\,$.

The case in which the monomial involves only the $\psi_i$'s is of course
covered by the Witten conjecture [Wi], proven by Kontsevich [Ko].
As for instance explained in [AC~2], this also allows to compute the
intersection numbers involving both $\ka1$ and the $\psi_i$'s.

It remains to compute the intersection numbers involving a boundary
class. Such a number may be thought of as the intersection of the remaining
classes on the corresponding boundary component. A problem with this
approach appears to be that most boundary components have singularities
that are not quotient singularities, which means that one cannot properly
do intersection theory on them. This problem is easily solved: we have
seen that each boundary component is the image under a finite map of
a moduli space of stable pointed curves or a product of two such spaces.
(The map almost always has degree $1$; the only exceptions are the degree $2$ 
map from $\mbar{g-1,n+2}$ to $\Dirr$ and, in the case $g$ even, the degree 
$2$ map from $\mbar{g/2,1}\times\mbar{g/2,1}$ to $\Delta_{g/2}\,$.)
So we wish to pull back the remaining divisor classes by means of this map.
If we can express the pull-backs in terms of the basic classes on the
new moduli space(s), we will be done, by induction on the dimension
of the moduli space.

So the whole point is to understand the pull-backs of the basic divisor
classes from $\mbar{g,n}$ to the moduli spaces occurring in its boundary
components. 

It is clear that the $\psi_i$'s pull back to $\psi_i$'s
on the new moduli space(s): to the first $n$ on $\mbar{g-1,n+2}\,$, and to the
$|N_1|$ $\psi_i$'s on $\mbar{g_1,|N_1|+1}$ and the $|N_2|$ $\psi_i$'s on 
$\mbar{g_2,|N_2|+1}$ that correspond to the points that are not identified
in the map to $\Delta_{g_i,N_i}\,$.

As explained in [AC~2], the class $\kappa_1$ pulls back to $\kappa_1$
on $\mbar{g-1,n+2}$ resp.~to the sum 
of the pull-backs of the $\kappa_1$'s from the two factors
on the product $\mbar{g_1,|N_1|+1}\times\mbar{g_2,|N_2|+1}\,$.

Pulling back a boundary divisor different from the one under consideration
to the new moduli space(s) is not difficult. The main point is that two
distinct boundary divisors intersect transversally in the universal
deformation space (see [DM]). It remains to identify the boundary divisors
on the new moduli space(s) that arise as the inverse image of the
intersection of the boundary divisor under consideration with a different one.

For example, in the case $\Dirr\,$, the pull-back of the class 
$\delta_{h,M}$ to $\mbar{g-1,n+2}$ is the sum 
$$\delta_{h-1,M\cup\{n+1,n+2\}} + \delta_{h,M}\,,$$
with some exceptions: when $n=0$ and $2h=g$, the two classes in the sum
are equal and the pull-back consists of that class just once;
when $n=0$ otherwise, $\delta_{h,\emptyset}$ 
has been denoted $\delta_{g-1-h,\{1,2\}}$
above; when $h=0$ or $h=g$, the first resp.~second summand is not
defined and should be omitted.

We now consider the pull-backs of boundary divisors to a product
$\mbar{g_1,|N_1|+1}\times\mbar{g_2,|N_2|+1}\,$.
The pull-back of $\dirr$ is the sum of the $\dirr$'s on the two factors. 
It remains to find the pull-backs of the boundary divisors parametrizing
reducible curves of a different type than the one under consideration.

We start with the case $n=0$. Hence $N_1=N_2=\emptyset$. We may assume
that $g_1\le g_2$ and we want to pull back the class $\delta_h\,$, with
$h\le g-h$ and $h\neq g_1\,$. The general curve in a component of the
intersection of the two boundary divisors is a chain consisting of
3 irreducible components. For the genera of the 3 components, there are
{\sl a priori\/} the following 4 possibilities:
\item{(1)} $[h,g_1-h,g_2]$ , occurring when $g_1>h$ ;
\item{(2)} $[h,g_2-h,g_1]$ , occurring when $g_2>h$ ;
\item{(3)} $[g_1,h-g_1,g-h]$ , occurring when $h>g_1$ ;
\item{(4)} $[g_2,h-g_2,g-h]$ , occurring when $h>g_2$ .

\no Here the second entry refers to the genus of the middle component;
note that $[a,b,c]$ and $[c,b,a]$ describe the same type of curves.

We observe that $h\le g/2\le g_2\,$. In fact $h<g_2\,$, since equality
implies $h=g_1\,$, which we have excluded. So (4) never occurs, while
(2) always occurs. We conclude that the pull-back of $\delta_h$ equals
$$\left\{
\eqalign{
&pr_1^*\delta_{g_1-h,\{1\}} + pr_2^*\delta_{g_2-h,\{1\}}
\qquad\hbox{in case }g_1>h\,;\cr
&pr_2^*\delta_{g_2-h,\{1\}} + pr_2^*\delta_{h-g_1,\{1\}}
\qquad\hbox{in case }h>g_1\,,\cr
}
\right.
$$
with one exception: when $h=g-h$ (implying $h>g_1$), the two summands
in the second line are equal, and the actual pull-back consists of that
class just once.

Now for the case $n>0$. We may assume that $1\in N_1$ and we want to pull
back the class $\delta_{h,M}\,$, different from $\delta_{g_1,N_1}\,$.
So $1\in M$ and $(h,M)\neq(g_1,N_1)$, but we no longer have $g_1\le g_2$
or $h\le g-h$. Again there are {\sl a priori\/} 4 possibilities for
the genera of the 3 components,
where as before the second entry refers to the genus of the middle component:
\item{(1)} $[h,g_1-h,g_2]$ , occurring when $g_1\ge h$ and $M\subset N_1\,$; 
\item{(2)} $[h,g_2-h,g_1]$ , occurring when $g_2\ge h$ and $M\subset N_2\,$;
\item{(3)} $[g_1,h-g_1,g-h]$ , occurring when $h\ge g_1$ 
and $N_1\subset M\,$;
\item{(4)} $[g_2,h-g_2,g-h]$ , occurring when $h\ge g_2$
and $N_2\subset M\,$.
\medskip
Observe that (2) never occurs, since $1\in M$ and $1\notin N_2\,$.
Note also that the other possibilities indeed yield stable curves in all
cases: in (1) and (3) the necessary condition $M\neq N_1$ when
$h=g_1$ is fulfilled, and in (4) the equality $M=N_2$ never occurs.
Finally, note that the types (1), (3) and (4) never coincide.

This means that the pull-back of $\delta_{h,M}\,$ consists of the sum
of 0, 1 or 2 of the classes from the following list, depending on
which conditions are satisfied:
\item{(1)} $pr_1^*\delta_{h,M}$ (pulled back from $\mbar{g_1,N_1\cup\{*\}}$)
when $g_1\ge h$ and $M\subset N_1\,$;
\item{(3)} $pr_2^*\delta_{h-g_1,M-N_1\cup\{1\}}$ (pulled back from
$\mbar{g_2,N_2\cup\{1\}}$) when $h\ge g_1$ and $N_1\subset M\,$;
\item{(4)} $pr_1^*\delta_{h-g_2,M-N_2\cup\{*\}}$ (pulled back from
$\mbar{g_1,N_1\cup\{*\}}$) when $h\ge g_2$ and $N_2\subset M\,$.
\medskip
\no Here we have identified the factors of the product of moduli spaces
by means of sets of marked points instead of just their number of
elements. In light of our convention to label a divisor parametrizing
reducible curves by the genus of the component containing 1 and by the
set of marked points on that component, it is natural to give the `extra'
point in case (3) the label 1, rather than $*$.

Finally, we have to deal with self-intersections of boundary divisors:
we need to pull back the class of a boundary divisor to the
corresponding (product of) moduli space(s). It is not difficult to deal
with this directly, but it is easier to use the fundamental identity
$$\ka1=12\la1-\delta+\psi$$
on $\mbar{g,n}\,$, where $\delta$ is the sum of the functorial classes
of the boundary divisors and $\psi$ is the sum of the $n$ $\psi_i$'s
(see [Co]). Namely, note that every divisor class on $\mbar{g,n}$
that we have discussed so far, occurs in this identity. So a given
boundary divisor class can be expressed as a linear combination of other
divisor classes, and if we know how to pull back the other classes, we will
also know how to deal with self-intersections. We have discussed the
pull-backs of $\ka1\,$, the $\psi_i$'s and the other boundary divisor
classes above, so we only need to determine the pull-back of $\la1\,$.
Note that the pull-back of the Hodge bundle to $\mbar{g-1,n+2}$ is an extension
of a trivial line bundle by the Hodge bundle in genus $g-1$, whereas the
pull-back of the Hodge bundle to $\mbar{g_1,|N_1|+1}\times\mbar{g_2,|N_2|+1}$
is the direct sum of the Hodge bundle in genus $g_1$ and the Hodge
bundle in genus $g_2\,$. (Cf.~[Kn~2] or [Co].)
Hence we find that the pull-back of $\la1$
to $\mbar{g-1,n+2}$ equals $\la1\,$, whereas the pull-back of $\la1$ to
$\mbar{g_1,|N_1|+1}\times\mbar{g_2,|N_2|+1}$ equals
$pr_1^*\la1+pr_2^*\la1\,$.

This determines then the pull-back of a boundary divisor class to its
corresponding (product of) moduli space(s). We find that 
in the case $n=0$ the pull-back of $\dirr$ to $\mbar{g-1,2}$ equals
$$-\psi_1-\psi_2+\dirr+\sum_{h=1}^{g-2}\delta_{h,\{1\}}\,,$$
whereas in the case $n>0$ the pull-back of $\dirr$ to $\mbar{g-1,n+2}$
equals
$$-\psi_{n+1}-\psi_{n+2}+\dirr+
\sum_{{0\le h\le g-1,\,1\in M\subset\nbar}
\atop{M\neq\nbar\,{\rm when}\,h=g-1}}
(\delta_{h,M\cup\{n+1\}}+\delta_{h,M\cup\{n+2\}})\,.
$$
In the case of a boundary divisor parametrizing reducible curves,
an actual self-intersection is much rarer. We find that the pull-back
of $\delta_{g_1,N_1}$ to 
$\mbar{g_1,N_1\cup\{*\}}\times \mbar{g_2,N_2\cup\{1\}}$ equals
$$
\left\{
\eqalign{
-pr_1^*\psi_{\{*\}}-pr_2^*\psi_1+pr_2^*\delta_{g_2-g_1\,,\{1\}}
&\qquad{\rm when}\quad n=0\quad{\rm and}\quad g_1<g_2\,;\cr
-pr_1^*\psi_{\{*\}}-pr_2^*\psi_1+pr_1^*\delta_{g_1-g_2\,,\,\nbar\cup\{*\}}
&\qquad{\rm when}\quad n>0,\quad N_1=\nbar\quad{\rm and}\quad g_1\ge g_2>0
\,;\cr
-pr_1^*\psi_{\{*\}}-pr_2^*\psi_1&\qquad{\rm otherwise}\,.\cr
}
\right.
$$
This finishes the theoretical description of the algorithm.
 
\bigskip
\no {\bf \S3. Implementation and results.}

\medskip
We have implemented the algorithm outlined in \S2 in 
Maple\footnote{$^1$}{Maple{\copyright} is a trademark
of the University of Waterloo and Waterloo Maple Software.}.
This first of all requires an implementation of the algorithm
for the computation of the intersection numbers of the $\psi_i$'s
on $\mbar{g,n}$ given in [Wi]. For this, we have gratefully used
the results of Chris Zaal 
who, using such an implementation, computed a table ([Za]) containing
all the intersection numbers
$$ 
\langle \tau_{d_1}\tau_{d_2}\cdots\tau_{d_n}\rangle
=\psi_1^{d_1}\psi_2^{d_2}\cdots\psi_n^{d_n}
$$
on $\mbar{g,n}$ with $g\le9$ for which all $d_i\ge2$ (hence
$n\le24$ since $\sum_{i=1}^n(d_i-1)=3g-3$).
The intersection numbers for which a $d_i$ equals 0 or 1 can be
computed from these by means of the string and dilaton equations.

{From} the intersection numbers of the $\psi_i$'s, one can determine
the intersection numbers of the classes $\ka{i}$ on $\mbar{g}$
introduced by Mumford in [Mu], as was briefly explained in [Wi].
Arbarello and Cornalba introduced in [AC~2] classes $\ka{i}$
on $\mbar{g,n}$ that generalize Mumford's classes, and they show
that the intersection numbers of the $\psi_i$'s determine the
intersection numbers of these $\ka{i}$'s as well as the `mixed'
intersection numbers of $\ka{i}$'s and $\psi_i$'s. So in particular
the intersection numbers of the divisor classes $\ka1$ and the $\psi_i$'s
are determined. It is easy to implement the calculation of these numbers
from the intersection numbers of the $\psi_i$'s (especially so
with a formula we learned from Dijkgraaf [Dij]),
and we have for instance calculated the numbers $\ka1^{3g-3}$ on $\mbar{g}$
for $g\le9$.

Naturally the various divisors on $\mbar{g,n}$ have to be ordered in some
consistent way. On $\mbar{g}$ we start with $\ka1\,$, followed by $\dirr$
and then the `reducible' boundary divisors, ordered by the minimum of the
genera of the 2 components, for a total of $[g/2]+2$ classes.
(The class $\la1$ was introduced only to deal with
self-intersections of boundary divisors and is not actually used in the 
program.)
When $n>0$, there are $(g+1)2^{n-1}+1$ classes: first
$\psi_1\,,\dots,\psi_n\,,\ka1$ and $\dirr\,$, then the reducible classes,
ordered first by the genus of the component that contains the point 1,
then by the number of points on that component, and finally by the
lexicographic ordering of subsets of $\nbar$ (of equal size and containing
1). (Recall that in genus 0
the class $\dirr$ is 0; it is included for convenience.)

In the case of a pull-back to a product of moduli spaces, it is
necessary to renumber the indices from $N_1\cup\{*\}$ as well as those
from $\{1\}\cup N_2\,$. For this we just use the natural ordering
of the elements of $\nbar\,$, taking $*$ as the $(n+1)$-st point.

The implementation of the actual algorithm is now rather straightforward.
Given a monomial in the divisor classes that contains at least one
boundary divisor, we order the classes as above, and pull back to the
(product of) moduli space(s) corresponding to the last occurring divisor.
Note that whenever the monomial contains several distinct boundary
divisors, we have a choice here; while the ordering chosen on $\mbar{g}$
is probably optimal, this is most likely not so when $n>0$. 
(Only recently have we experimented with another ordering, which indeed
appeared to be better; the idea is that in the reducible case one should
try to pull back to two moduli spaces parametrizing curves of
approximately equal genus and number of points.)

In case the last divisor is $\Dirr\,$, we find a homogeneous polynomial
of degree $3g-4+n$ in the divisor classes on $\mbar{g-1,n+2}\,$.
After expanding, it is a sum of monomials (with coefficients) in the
divisor classes; these are evaluated by means of the (heavily
recursive) algorithm.

In the reducible case, we also find a homogeneous polynomial of
degree $3g-4+n$, but this time in two sets of variables, the divisors
on $\mbar{g_1,|N_1|+1}$ and those on $\mbar{g_2,|N_2|+1}\,$. Many
of the monomials will not have the correct degree $3g_1-2+|N_1|$ in the
first set of variables and are 0 for trivial reasons. The others
automatically have degree $3g_2-2+|N_2|$ in the second set of variables;
writing such a monomial as $c\cdot M_1\cdot M_2\,$, where $c$ is the
coefficient of the monomial and $M_i$ is the
monic monomial in the $i$-th set of variables,
it contributes $c\cdot a(M_1)\cdot a(M_2)$, where $a(M_i)$ is the
result of applying the algorithm to $M_i$ on $\mbar{g_i,|N_i|+1}\,$.

Using this implementation of the algorithm, we have computed e.g.~the
28 intersection numbers of $\ka1$, $\dirr$ and $\delta_1$ on $\mbar3\,$,
confirming the results of [Fa~1]. However, because of its heavily recursive 
character, the algorithm becomes impracticable already in the computation
of certain intersection numbers on $\mbar4\,$. Most intersection numbers
are still easy to compute, but especially the numbers $\ka1^{9-i}\dirr^i$ with 
$i$ large take a long time. It is quite clear why: firstly, a pull-back
to $\mbar{g-1,n+2}$ is the `worst case', since the changes in genus and in
dimension of the moduli space are minimal, while the number of points
increases by 2, so that the number of divisors increases by a factor of 
almost 4;
secondly, as we saw in \S2, the pull-back of $\dirr$ to
$\mbar{g-1,n+2}$ involves by far the highest number of terms.

Observe however that the class $\dirr$ is a pull-back from $\mbar{g}$
resp.~$\mbar{1,1}\,$. This first of all implies
$$\dirr^{m+1}=0,
$$
where $m=\max(g,3g-3)$, for all $g\ge0$. A systematic use of this 
identity already saves considerable time.
Moreover, any product involving only $\ka1\,$, the $\psi_i$'s and 
$\dirr$ can be pushed down to $\mbar{g}$ resp.~$\mbar{1,1}\,$, with the use
of the projection formula and the formulas in [AC~2]. 
This leads to intersection numbers on those spaces of monomials in
$\dirr$ and the higher $\ka{i}$ mentioned before.
As observed by Arbarello and Cornalba, the(ir) $\ka{i}$ behave very well
under pull-back, and it is clear that all intersection numbers 
involving divisors as well as the $\ka{i}$ can be computed by means
of an algorithm almost identical to the one in \S2. Even the
implementation is easy to adapt. The point is that this greatly
simplifies the calculation of 
the numbers involving only $\ka1\,$, the $\psi_i$'s and
$\dirr\,$: in the naive implementation,
the complexity of a calculation increases at least exponentially with the
number of points, and one would have to calculate certain numbers on
$\mbar{0,2g}$ in order to get all numbers on $\mbar{g}\,$; but now many
of the hardest numbers can be computed using at most $2$-pointed curves
at all stages of the computation.

With these simple changes implemented, the calculation of many more
numbers becomes practical. We have calculated all intersection numbers
of divisors on $\mbar{g}$ for $g\le6$ as well as on $\mbar{3,1}$ and
$\mbar{4,1}\,$. To obtain the 2 numbers 
$\ka1\dirr^{14}$ and $\dirr^{15}$ on $\mbar{6}\,$, we used
the relations $\la1^{13}\dirr\ka1=0$ and $\la1^{13}\dirr^2=0$,
consequences of the geometrically obvious relation
$\la1^{13}\dirr=0$.

This section of the paper would hardly be complete without some actual
intersection {\sl numbers\/}. Here are a few:
$$\displaylines{
{\rm On}\,\mbar4\,:\qquad\dirr^9={{-251987683}\over{4320}}\,,\quad
\la1^9={1\over{113400}}\,;\cr
{\rm On}\,\mbar5\,:\qquad\dirr^{12}={{-1766321028967}\over{6048}}\,,\quad
\la1^{12}={{31}\over{680400}}\,;\cr
{\rm On}\,\mbar6\,:\qquad\dirr^{15}={{-32467988437272065977}
\over{7257600}}\,,\quad
\la1^{15}={{431}\over{481140}}\,.\cr
}
$$
We computed the number $\la1^9$ on $\mbar4$ in [Fa~2] by a completely
different (ad hoc) method. Calculating it with the algorithm amounts
to calculating all $220$ intersection numbers of divisors on $\mbar4\,$,
so this provides a nice check of the implementation.
 
\vfill\eject
%\bigskip
\no {\bf \S4. The class of the locus of Jacobians and other applications.}

\medskip
The calculation of $\la1^9$ on $\mbar4$ in [Fa~2] was used there to
obtain the well-known result that the class of the locus $\j4$ of
Jacobians of curves of genus 4 in the moduli space $\a4$ of
principally polarized abelian varieties of dimension 4 equals
$8\la1\,$. Using the computation of $\la1^{12}$ on $\mbar5$ as well as
some computations in the tautological ring of $\m5\,$ as in [Fa~3],
we could determine the class of $\j5$ in $\a5\,$, as we will explain
in a moment. Recently we realized that Mumford's formula [Mu] for the
Chern character of the Hodge bundle on $\mbar{g}\,$, together with
an algorithm similar to the one discussed in \S2, enable one to
compute the class of $\j{g}$ in $\a{g}$ for {\sl all\/} $g$, at least
in principle. In practice, we have carried this out for $g\le7$. These
results will be discussed here as well.

First we recall the set-up and explain what we
mean by ``the class of the locus of Jacobians in the moduli
space of principally polarized abelian varieties": this will be the
class in the {\sl tautological ring\/} of $\tilde\a{g}\,$, the
$\Q$-subalgebra of the cohomology ring of a toroidal compactification
$\tilde\a{g}$ of $\a{g}$ generated by the Chern classes $\la{i}$ of the
Hodge bundle $\E$ on $\tilde\a{g}\,$. From [Mu], \S5 we know that the relation
$$
(1-\la1+\la2-\la3+\dots+(-1)^g\la{g})\,(1+\la1+\la2+\la3+\dots+\la{g})=1
$$
holds; equivalently, $\ch{2k}(\E)=0$ for all $k\ge1$.
The tautological ring is in fact the quotient of $\Q[\la1,\dots,\la{g}]$
by the ideal generated by the homogeneous components of the relation
above and is thus a complete intersection ring. A detailed description
of the tautological ring may be found in [vdG]. In particular, the relation
with the cohomology ring of the {\sl compact dual\/} of the Siegel upper
half space via Hirzebruch's proportionality principle is explained there.
This includes the fundamental identity
$$
\prod_{i=1}^g\,\la{i}=\prod_{i=1}^g\,{{|B_{2i}|}\over{4i}}
$$
that enables one to compute intersection numbers in the tautological ring
of $\tilde\a{g}\,$.

Denoting by $t:\mbar{g}\to\tilde\a{g}$ the extended Torelli morphism
and its image by $\tilde\j{g}\,$, we are after the functorial class
$[\tilde\j{g}]_Q$ of the locus of (generalized) Jacobians, which is
one half its usual fundamental class. In other words, we wish to
determine ${1\over2}t_*1$, since a generic curve of genus at least 3
has no non-trivial automorphisms, while the generic p.p.a.v.~and the
generic Jacobian of dimension at least 3 have two automorphisms.

It is important to point out that this is {\sl not\/} what we
actually compute. We do not know whether the class of the locus
of Jacobians actually lies in the tautological ring (our feeling
is that this is not the case for $g$ large enough, but we cannot even
think of a method to decide this). Instead, we compute the
{\sl projection\/} of this class in the tautological ring; this is 
well-defined by the perfect pairing in the tautological ring and the
cohomology ring of $\tilde\a{g}\,$. (In other words, we compute the
class modulo a class $X$ that pairs zero 
with all classes in the tautological ring of the complementary dimension.)

The method of computation of {\sl this\/} class is the following.
It is a class in the tautological ring of $\tilde{\a{g}}$ of dimension
$3g-3$, hence of codimension $c={{g-2}\choose2}$. Write it as a linear
combination with unknown coefficients $a_i$ of the elements $s_i$
of a basis of the degree-$c$ part of the tautological ring:
$$
[\tilde{\j{g}}]_Q={1\over2}t_*1=a_1s_1+a_2s_2+\dots+a_ks_k\,(+X).
$$
(A natural basis is for instance the collection of square-free
monomials of degree $c$ in $\la1\,,\dots,\la{g}\,$.) To compute
the coefficients $a_i\,$, we have to evaluate the $k$ monomials
$\Lambda_i$ of a basis of the degree-$(3g-3)$ part of the tautological
ring on this class. The evaluation of expressions $\Lambda_is_j$ in the
tautological ring of $\tilde{\a{g}}$ uses the relations between the
$\la{i}$ and the proportionality relation stated above.
So the expressions $\Lambda_i(a_1s_1+\dots+a_ks_k)$ yield $k$ linearly
independent rational linear combinations of the unknowns $a_i\,$.
The values of these expressions can be determined as 
${1\over2}t_*(t^*\Lambda_i)=({1\over2}t_*1)\cdot\Lambda_i\,$,
if we know how to evaluate $t^*\Lambda_i$ on $\mbar{g}\,$.

The simplest non-trivial example is $g=4$. Here $c=1$, a basis in
codimension 1 is $\la1\,$, a basis in codimension 9 is $\la1^9$
(or any non-zero monomial in the $\la{i}$), so to compute the
class of $\tilde\j4$ in $\tilde\a4$ we only need $\la1^9$ on
$\mbar4\,$. We have seen above that this can be evaluated e.g.~by
means of the implementation of the algorithm; we find the
well-known result that $[\tilde\j4]_Q=8\la1\,$.

The situation for $g=5$ is more interesting. Here $c=3$ and a basis
in codimension 3 is given by $\la1\la2\,(={1\over2}\la1^3)$ and $\la3\,$.
We need to evaluate 2 independent monomials of degree 12 in the $\la{i}$
on $\mbar5\,$. The algorithm will naturally yield only the number
$\la1^{12}$ (whose value we gave at the end of \S3). However, we can use
the simple observation that the 
class $\la{g}\la{g-1}$ vanishes
on the boundary $\mbar{g}-\m{g}$ (see
[Fa~3] or [Fa~4]). As a corollary, the numbers $\la1\la2\la4\la5$
and $\la3\la4\la5$ on $\mbar5$ satisfy the same relation as the
classes $\la1\la2$ and $\la3$ in the 1-dimensional degree-3 part
of the {\sl tautological ring\/} $R^*(\m5)$ of $\m5\,$. This relation
was worked out in [Fa~3]: 
$10\la3=3\la1\la2\,$. A quick calculation in the tautological ring
of $\tilde\a5$ shows that this implies that the class of the
locus of Jacobians satisfies
$$
[\tilde\j5]_Q=a\,(3\la1\la2-2\la3).
$$
Another such calculation, using $\la1^{12}={{31}\over{680400}}$ on
$\mbar5\,$, shows then that $a=24$, hence
$$
[\tilde\j5]_Q=72\la1\la2-48\la3\,.
$$

For higher genus, we find with this method only some of the
coefficients, not all of them. We start with a general formula:

\proclaim Conjecture 1. In the basis of monic square-free monomials
in $\la1\,,\dots,\la{g}$ of degree ${{g-2}\choose2}$, the coefficient
$C_{1,2,\dots,g-3}$
of $\la1\la2\cdots\la{g-3}$ in the (projected) class $[\tilde\j{g}]_Q$
equals
$${1\over{2g-2}}\,\prod_{i=1}^{g-2}\,{2\over{(2i+1)|B_{2i}|}}\,. $$
So it equals 1, 8, 72, 384, 768 resp.~for $g=3,\dots,7$, 
while it is not an integer for any larger value of $g$.

\no The conjecture is true for $g\le15$, as it is a consequence
of a conjectural formula in [Fa~3] for the number $\la{g-1}^3\,$
on $\mbar{g}$ that is proven for $g\le15$.
To derive it from that formula, 
simply note that $\la{g-1}^3=2\la{g-2}\la{g-1}\la{g}\,$ and that
$\la1\la2\cdots\la{g-3}$ is the unique monomial in the basis that
pairs non-zero with $\la{g-2}\la{g-1}\la{g}\,$.

Note that we could have used this formula instead of that for 
$\la1^{12}\,$ to compute the class of $[\tilde\j5]_Q\,$.

For $g=6$
we find the coefficient of $\la2\la4$ using the
relation between $\la1\la3$ and $\la4$ in $R^*(\m6)$:
$$C_{2,4}=-3C_{1,2,3}=-1152.
$$
The knowledge of $\la1^{15}$ (see \S3) gives a non-trivial relation between 
the remaining coefficients of $\la1\la5$ and $\la6$:
$$
C_6+16C_{1,5}={{7336704}\over{691}}\,.
$$

A new ingredient will be required to solve for these coefficients.
It is provided by Mumford's formula [Mu] for the Chern character of the
Hodge bundle on $\mbar{g}\,$:
$$
ch(\E)=g+\sum_{i=1}^{\infty}{{B_{2i}}\over{(2i)!}}
\left[\ka{2i-1}+{1\over2}\sum_{h=0}^{g-1}i_{h,*}\left(
K_1^{2i-2}-K_1^{2i-3}K_2+\dots+K_2^{2i-2}\right)\right]\,.
$$
(Note that we use a different convention for the Bernoulli numbers:
$B_2={1\over6}$, $B_4={{-1}\over{30}}$ etc.)
Here $i_0:\mbar{g-1,2}\to\Dirr\subset\mbar{g}$ and 
$i_h:\mbar{h,1}\times\mbar{g-h,1}\to\Delta_h\subset\mbar{g}$ are the
natural maps, and $K_i$ is the first Chern class of the relative
cotangent line bundle at the $i$-th point.

The formula is ideally suited for a recursive computation of the
intersection numbers of the $\ka{i}$ and the $\ch{2j-1}(\E)$.
Namely, suppose given a monomial in those classes of degree
$3g-3+n$ on $\mbar{g,n}\,$. If only $\ka{i}$ occur, we can proceed
as explained in \S2, using the Witten conjecture (Kontsevich's theorem).
In any case, as the Hodge bundle is a pull-back from $\mbar{g}$
resp.~$\mbar{1,1}\,$, we can push down the expression to $\mbar{g}$
resp.~$\mbar{1,1}$ and obtain a sum of similar expressions.
In genus 1, we only need the well-known equalities
$$
\ch1(\E)=\la1=\ka1=\psi_1=\ll\tau_1\rr={1\over{24}}\,.
$$
So assume the genus is at least 2 and the monomial on $\mbar{g}$ contains
at least one $\ch{2j-1}(\E)$. Take the highest odd Chern character
component that occurs, and expand it using Mumford's formula.
In the first term, that $\ch{2k-1}$ is replaced 
by a $\ka{2k-1}$ (up to a factor), so it is determined inductively.
The other terms involve expressions in the $K_i$ pushed forward
via the maps $i_h\,$. The point is that these can be written as 
push-forwards from $\mbar{g-1,2}$ resp.~$\mbar{h,1}\times\mbar{g-h,1}$
of intersection numbers of the classes $K_i\,$, $\ka{j}$ and
Chern character components of the Hodge bundles in genus $g-1$ resp.~genera
$h$ and $g-h$. This is clear from the Arbarello-Cornalba formulas
for $i_h^*\ka{j}$ and 
the fact that $i_0^*(\E)$ is the extension of a trivial line bundle
by the Hodge bundle in genus $g-1$, while for $h$ positive 
$i_h^*(\E)$ is the direct sum of the Hodge bundles in genera $h$ and $g-h$.

After expanding and omitting the terms that are 0 for dimension reasons,
we find an expression in intersection numbers of the classes just
mentioned on spaces of 1- or 2-pointed curves. These can be pushed down
again to intersection numbers of $\ka{i}$ and $\ch{2j-1}(\E)$
on $\mbar{h}$ (with $h<g$) or $\mbar{1,1}\,$. By induction on the
genus, these numbers are known.

Having discussed the implementation of the divisor algorithm in some
detail, we content ourselves with saying that the implementation of the
new algorithm proceeds along similar lines and is considerably easier.

Because the $\la{i}$ can be expressed in the $\ch{j}(\E)$, this means
that all the intersection numbers of the $\la{i}$ on $\mbar{g}$ can be
computed recursively. By the discussion above, it follows that
for {\sl all\/} $g$ the projection in the tautological ring
of $\tilde\a{g}$ of the class of the locus of Jacobians can be
computed, at least in principle.

Currently, we have carried this out for $g\le7$. For $g=6$, only
one more relation was required. Either one of the two relations
following from
$$
\la2\la3\la4\la6={{1697}\over{2988969984000}}\,,\qquad
\la1\la2\la3\la4\la5={{150719}\over{15692092416000}}\,
$$
suffices to solve for $C_6$ and $C_{1,5}$ (the first relation involves
$C_{1,5}$ only). The result is
$$\eqalign{
[\tilde\j6]_Q&=384\la1\la2\la3-1152\la2\la4+{{474048}\over{691}}\la1\la5
-{{248064}\over{691}}\la6  	\cr
{}&=2^73\left(\la1\la2\la3-3\la2\la4+{{2469}\over{1382}}\la1\la5
-{{646}\over{691}}\la6\right). 	\cr
}
$$
It may be worthwhile to point out the relation $15C_6+28C_{1,5}=2^93^3$.

In genus 7 the result is:
$$\eqalign{
[\tilde\j7]_Q&=768\la1\la2\la3\la4-6912\la2\la3\la5+
{{2209152}\over{691}}\la1\la4\la5
+{{7522176}\over{691}}\la1\la3\la6   \cr
&\qquad{}-{{8842752}\over{691}}\la4\la6
+{{968832}\over{691}}\la3\la7-{{3276672}\over{691}}\la1\la2\la7\,.\cr
}
$$
As stated, this result is not very pretty; perhaps the class looks better
in a different basis. We would like to point out though that the class
can be computed with any choice of 7 independent monomials of degree 18
in the $\la{i}\,$; once this is accomplished, the values of all other
such monomials can be determined by means of an easy calculation in
the tautological ring of $\tilde\a7\,$. In particular, one can choose
`easy' monomials to compute the class, and get the `hard' ones for free.
In this way we computed for instance $\la1^{18}$ on $\mbar7\,$:
$$
\la1^{18}={{32017001}\over{638512875}}\,.
$$

Finally, we have another general formula:
\proclaim Conjecture 2. In the basis of monic square-free monomials
in $\la1\,,\dots,\la{g}$ of degree ${{g-2}\choose2}$, the coefficient
$C_{2,3,\dots,g-4,g-2}$
of $\la2\la3\cdots\la{g-4}\la{g-2}$ in the (projected) class $[\tilde\j{g}]_Q$
equals
$$\left({{g(2g-2)}\over{12}}-2^{g-3}\right)\,C_{1,2,\dots,g-3}
={g\over{12}}\,\prod_{i=1}^{g-2}\,{2\over{(2i+1)|B_{2i}|}}
-{1\over{4g-4}}\,\prod_{i=1}^{g-2}\,{4\over{(2i+1)|B_{2i}|}}\,.$$
So it equals -48, -1152, -6912 resp.~for $g=5,6,7$,
while it is not an integer for any larger value of $g$.

\no Again, the conjecture is true for $g\le15$; 
it is a consequence of Conjecture 1 and
of the conjectural formula $\ka1\la{g-3}=g(2g-2)\la{g-2}$ in
$R^{g-2}(\m{g})$ that can be proven for $g\le15$ using the results of [Fa~3]
(note that the only 2 monomials in the basis that pair non-zero with
$\la1\la{g-3}\la{g-1}\la{g}$ are $\la2\la3\cdots\la{g-4}\la{g-2}$
and $\la1\la2\dots\la{g-3}$).

\bigskip
\no {\bf Acknowledgements and final remarks.}

\medskip
The research in the first 3 sections and the determination of the
classes of $\tilde\j4$ and $\tilde\j5$ was carried out at the
Universiteit van Amsterdam and was made possible
by a fellowship of the Royal Netherlands Academy of Arts and Sciences.

Rahul Pandharipande's remark that led directly to the algorithm in \S2
was made during a visit of the author to the University of Chicago
in April/May 1995. The author would like to thank Bill Fulton for the
invitation to visit and the University of Chicago for partial support.
Both this remark and an analogous remark by Eduard Looijenga half a year
earlier (that I didn't get around to working out during that half year,
unfortunately) were sparked by talks in which I discussed the problem
of determining the intersection numbers of divisors on $\mbar{g}$
and the partial results I had obtained (cf.~[Fa~2], \S5).
I am most grateful to them for these remarks. I also want to thank
Chris Zaal whose results I used heavily in my computations.

The research in the remaining part of \S4 was carried out at the
Institut Mittag-Leffler of the Royal Swedish Academy of Sciences
during the academic year 1996/97
devoted to {\sl Enumerative geometry and its interaction
with theoretical physics\/}. 
My stay at the Mittag-Leffler Institute was supported by a grant from
the G\"oran Gustafsson Foundation for Mathematics, Physics and
Medicine to Torsten Ekedahl. I would like to thank the organizers
of the program and the staff of the institute for making it
a wonderful year.

Finally, some remarks about possible uses of the implemented algorithms.
The algorithm for computing intersection numbers involving Chern classes
of the Hodge bundle (described in \S4) has already found applications
in the recent work of Graber and Pandharipande [GP] and
Kontsevich and Manin [KM].
For moduli spaces $\mbar{g,n}$ of small dimension and with not too many
divisors, one can use the divisor algorithm to
determine the part of the cohomology ring generated by divisors.
This makes e.g.~an intersection calculation on $\mbar{1,4}$ as in [Ge]
relatively painless. As in that paper, there may be applications
to computing Gromov-Witten invariants.
In fact, the algorithms involving the $\ka{i}$ and the $\ch{2j-1}$
allow to do calculations that include those classes; this covers
e.g.~the Chow ring of $\mbar3\,$. It is even possible to include
arbitrary boundary strata as {\sl module\/} generators, by pulling back
all other classes to the corresponding product of moduli spaces via
a sequence of maps, each identifying a single pair of points.
Writing algorithms that can handle {\sl intersections\/} of arbitrary boundary
strata will be considerably more difficult, however.

The various algorithms, and some tables of intersection numbers computed
with it, are available from the author by e-mail, although only 
the initial version of the divisor algorithm is currently available in
a user-friendly format.

\vfill\eject
%\bigskip
\noindent {\bf References.\/}
\medskip
\pn [AC 1] E.~Arbarello, M.~Cornalba, {\it The Picard groups of the moduli
spaces of curves\/}, Topology 26 (1987), 153--171.
\pn [AC 2] E.~Arbarello, M.~Cornalba, 
{\sl Combinatorial and algebro-geometric cohomology classes on the moduli
spaces of curves\/}, J.~Alg.~Geom. 5 (1996), 705-749.
\pn [Co] M.~Cornalba, {\it On the projectivity of the moduli spaces of 
curves\/}, J.~Reine Angew.~Math. 443 (1993), 11--20.
\pn [DM] P.~Deligne, D.~Mumford, {\it The irreducibility of the space of 
curves of given genus\/}, Publ.~Math.~IHES 36 (1969), 75--109.
\pn [Dij] R.~Dijkgraaf, {\sl Some facts about tautological classes\/},
private communication (November 1993).
\pn [Fa 1] C.~Faber, {\sl Chow rings of moduli spaces of curves I: The Chow
ring of $\mbar3$\/}, Ann.~of Math. 132 (1990), 331--419.
\pn [Fa 2] C.~Faber, {\sl Intersection-theoretical computations on
$\mbar{g}\,$\/}, in {\sl Parameter Spaces\/} (Editor P.~Pragacz),
Banach Center Publications, Volume 36, Warszawa 1996.
\pn [Fa 3] C.~Faber, {\sl A conjectural description of the tautological 
ring of the moduli space of curves\/}, preprint 1996
(available from http://www.math.okstate.edu/preprint/1997.html).
\pn [Fa 4] C.~Faber, {\sl A non-vanishing result for the 
tautological ring of $\Mg\,$\/}, preprint 1995
(available from http://www.math.okstate.edu/preprint/1997.html).
\pn [vdG] G.~van der Geer, {\sl Cycles on the moduli space of abelian
varieties\/}, alg-geom/9605011.
\pn [Ge] E.~Getzler, {\sl Intersection theory on $\mbar{1,4}$ and
elliptic Gromov-Witten invariants\/}, alg-geom/9612004.
\pn [GP] T.~Graber, R.~Pandharipande, {\sl Localization of virtual
classes\/}, alg-geom/9708001.
\pn [Ha 1] J.~Harer, {\sl The second homology group of the mapping class
group of an orientable surface\/}, Invent.~Math. 72 (1983), 221--239.
\pn [Ha 2] J.~Harer, {\sl The cohomology of the moduli space of curves\/},
in {\sl Theory of moduli\/}, 138--221, LNM 1337, Springer, Berlin-New York 1988.
\pn [Ke] S.~Keel, {\sl Intersection theory of moduli space of stable
$N$-pointed curves of genus zero\/}, Trans.~A.M.S. 330 (1992), 545--574.
\pn [Kn 1] F.~Knudsen, {\sl The projectivity of the moduli spaces of
stable curves, II: The stacks $\m{g,n}\,$\/}, Math.~Scand. 52 (1983), 161--199.
\pn [Kn 2] F.~Knudsen, {\sl The projectivity of the moduli spaces of
stable curves, III: The line bundles on $\m{g,n}$ and a proof of the
projectivity of ${\overline{M}}_{g,n}$ in characteristic 0\/},
Math.~Scand. 52 (1983), 200--212.
\pn [Ko] M.~Kontsevich, {\sl Intersection theory on the moduli space
of curves and the matrix Airy function\/}, Comm.~Math.~Phys. 147 (1992),
1--23.
\pn [KM] M.~Kontsevich, Yu.~I.~Manin, {\sl Relations between the correlators
of the topological sigma-model coupled to gravity\/}, alg-geom/9708024.
\pn [Mu] D.~Mumford, {\sl Towards an enumerative geometry of the moduli
space of curves\/}, in {\sl Arithmetic and Geometry\/} (Editors M.~Artin
and J.~Tate), Part II, Progress in Math., Volume 36,
Birkh\"auser, Basel 1983.
\pn [Wi] E.~Witten, {\sl Two dimensional gravity and intersection theory
on moduli space\/}, Surveys in Diff.~Geom. 1 (1991), 243--310.
\pn [Za] C.~Zaal, Maple procedures for computing the
intersection numbers of the $\psi_i$'s and a table of the numbers for
$g\le9$. Universiteit van Amsterdam 1992. (Available by e-mail from
zaal@wins.uva.nl or from the author.)
 
\bigskip
\no Department of Mathematics, Oklahoma State University, Stillwater,
Oklahoma 74078-1058, U.S.A.
 
\medskip\no
e-mail: cffaber@math.okstate.edu
 
\bye